\documentclass[10pt]{amsart}
\input{epsf}

\usepackage{amssymb,latexsym}
\topskip  \numberwithin{equation}{section}
\newtheorem{theorem}{Theorem}[section]
\newtheorem{proposition}[theorem]{Proposition}
\newtheorem{corollary}[theorem]{Corollary}

\newtheorem{example}[theorem]{Example}

\begin{document}

\pagenumbering{arabic}
\pagestyle{headings}
\def\sof{\hfill\rule{2mm}{2mm}}
\def\ls{\leq}
\def\gs{\geq}
\def\qq{{\bold q}}
\def\txx{{\frac1{2\sqrt{x}}}}
\def\tx{{\left(\txx\right)}}
\def\ttx{\left({{\frac{1}{2x}}}\right)}
\def\Bn{\mathcal{P}_n}
\def\mn{\mbox{-}}
\def\vn{\varnothing}
\def\vr{\varnothing}
\def\mnk{\mn k}
\def\pk{{1\mn2\mn\cdots\mnk}}
\def\pkk{{1\mn2\mn\cdots\mn(k-1)}}
\def\SS{\frak S}

\title{Restricted $132$-alternating permutations and Chebyshev polynomials}

\maketitle

\begin{center}Toufik Mansour \footnote{Research financed by EC's
IHRP Programme, within the Research Training Network "Algebraic
Combinatorics in Europe", grant HPRN-CT-2001-00272}
\end{center}

\begin{center}{Department of Mathematics, Chalmers University of
Technology, S-41296 G\"oteborg, Sweden

        {\tt toufik@math.chalmers.se} }
\end{center}
\section*{Abstract}
A permutation is said to be {\em alternating} if it starts with
rise and then descents and rises come in turn. In this paper we
study the generating function for the number of alternating
permutations on $n$ letters that avoid or contain exactly once
$132$ and also avoid or contain exactly once an arbitrary pattern
on $k$ letters. In several interesting cases the generating
function depends only on $k$ and is expressed via Chebyshev
polynomials of the second kind. \thispagestyle{empty}
\section{Introduction}
The main goal of this paper is to give analogues of enumerative
results on certain classes of permutations characterized by
pattern-avoidance in the symmetric group $\SS_n$. In the set of
alternating permutations we identify classes of restricted
alternating permutations with enumerative properties analogous to
results on permutations. More precisely, we study generating
functions for the number of alternating permutations that avoid or
contain exactly once $132$ and also avoid or contain exactly once
an arbitrary permutation $\tau\in\SS_k$. In the remainder of this
section, we present a brief account of earlier works which
motivated our investigation, we give the basic definitions used
throughout the paper, and we summarize our main results.

\subsection{Classical patterns} Let $\alpha\in\SS_n$ and $\tau\in
\SS_k$ be two permutations. We say that $\alpha$ {\it contains\/}
$\tau$ if there exists a subsequence $1\ls i_1<i_2<\dots<i_k\ls n$
such that $(\alpha_{i_1}, \dots,\alpha_{i_k})$ is order-isomorphic
to $\tau$; in such a context $\tau$ is usually called a {\it
pattern\/}. We say that $\alpha$ {\it avoids\/} $\tau$, or is
$\tau$-{\it avoiding\/}, if such a subsequence does not exist. The
set of all $\tau$-avoiding permutations in $\SS_n$ is denoted
$\SS_n(\tau)$. For an arbitrary finite collection of patterns $T$,
we say that $\alpha$ avoids $T$ if $\alpha$ avoids any $\tau\in
T$; the corresponding subset of $\SS_n$ is denoted $\SS_n(T)$.

While the case of permutations avoiding a single pattern has
attracted much attention, the case of multiple pattern avoidance
remains less investigated. In particular, it is natural, as the
next step, to consider permutations avoiding pairs of patterns
$\tau_1$, $\tau_2$. This problem was solved completely for
$\tau_1,\tau_2\in \SS_3$ (see \cite{SS}), for $\tau_1\in \SS_3$
and $\tau_2\in \SS_4$ (see \cite{W}), and for $\tau_1,\tau_2\in
\SS_4$ (see \cite{B1,Km} and references therein). Several recent
papers \cite{CW,MV1,Kr,MV2,MV3,MV4} deal with the case $\tau_1\in
\SS_3$, $\tau_2\in \SS_k$ for various pairs $\tau_1,\tau_2$.
Another natural question is to study permutations avoiding
$\tau_1$ and containing $\tau_2$ exactly $t$ times. Such a problem
for certain $\tau_1,\tau_2\in \SS_3$ and $t=1$ was investigated in
\cite{R}, and for certain $\tau_1\in \SS_3$, $\tau_2\in \SS_k$ in
\cite{RWZ,MV1,Kr}. The tools involved in these papers include
Fibonacci numbers, Catalan numbers, Chebyshev polynomials,
continued fractions, and Dyck words.

{\it Fibonacci numbers} are defined by $F_0=0$, $F_1=1$, and
$F_n=F_{n-1}+F_{n-2}$ for all $n\geq2$.

{\it Catalan numbers} are defined by
$C_n=\frac{1}{n+1}\binom{2n}{n}$ for all $n\geq0$. The generating
function for the Catalan numbers is given by
$C(x)=\frac{1-\sqrt{1-4x}}{2x}$.

{\it Chebyshev polynomials of the second kind\/} (in what follows
just Chebyshev polynomials) are defined by
$U_r(\cos\theta)=\frac{\sin(r+1)\theta}{\sin\theta}$ for $r\geq0$.
Evidently, $U_r(t)$ is a polynomial of degree $r$ in $t$ with
integer coefficients, and satisfy
\begin{equation}
U_0(t)=1,\ U_1(t)=2t,\ \mbox{and}\ U_r(t)=2tU_{r-1}(t)-U_{r-2}(t)\
\mbox{for all}\ r\geq2.\label{reccheb}
\end{equation}
Chebyshev polynomials were invented for the needs of approximation
theory, but are also widely used in various other branches of
mathematics, including algebra, combinatorics, and number theory
(see \cite{Ri}). Apparently, for the first time the relation
between restricted permutations and Chebyshev polynomials was
discovered by Chow and West in \cite{CW}, and later by Mansour and
Vainshtein \cite{MV1,MV2,MV3,MV4}, Krattenthaler \cite{Kr}. These
results related to a rational function
\begin{equation}
R_k(x)=\frac{U_{k-1}\tx}{\sqrt{x}U_k\tx}, \label{rcheb}
\end{equation}
for all $k\geq 1$. It is easy to see that for any $k$, $R_k(x)$ is
rational in $x$ and satisfies
\begin{equation}R_k(x)=\frac{1}{1-xR_{k-1}(x)},\label{rk}\end{equation}
for all $k\geq1$.

\subsection{Generalized patterns} In \cite{BS}, there were introduced
generalized permutation patterns that allow the requirement that
two adjacent letters in a pattern must be adjacent in the
permutation. We write a classical pattern with dashes between any
two adjacent letters of the pattern, say $1342$, as
$1\mn3\mn4\mn2$, and if we write, say $24\mn3\mn1$, then we mean
that if this pattern occurs in a permutation $\pi$, then the
letters in $\pi$ that correspond to $2$ and $4$ are adjacent (for
more details, see \cite{C}). For example, the permutation
$\pi=35421$ has only two occurrences of the pattern $23\mn1$,
namely the subsequences $352$ and $351$, whereas $\pi$ has four
occurrences of the pattern $2\mn3\mn1$, namely the subsequences
$352$, $351$, $342$ and $341$. In \cite{MV1,Mg1,Mg2,Mg3}, there
were studied the number of permutations that avoid or contain
exactly once $1\mn3\mn2$ and that also avoid or contain exactly
once an arbitrary generalized pattern. The tools involved in these
papers include Fibonacci numbers, Catalan numbers, Chebyshev
polynomials, and continued fractions.

\subsection{Alternating permutations} A permutation
is said to be alternating if it starts with \emph{rise} and then
\emph{descents} and rises come in turn. In other words, an
alternating permutation $\pi\in \SS_n$ satisfies
$\pi_{2j-1}<\pi_{2j}>\pi_{2j+1}$ for all $1\leq j\leq [n/2]$, that
is to say its \emph{rise} (resp. \emph{descent}) is equal to an
odd (resp. even) index. We denote the set of all alternating
permutations on $n$ letters by $A_n$. Other names that authors
have used for these permutations are \emph{zig-sag permutation}
and \emph{up down permutations}. The determination of the number
of alternating permutations for the set $\{1,2,\ldots,n\}$ (or on
$n$ letters) is known as Andr\'e's problem (see~\cite{AP1, AP2}).
An example of an alternating permutation is $14253$. The number of
alternating permutations on $n$ letters, for $n=1,2,\ldots,10$, is
$1$, $1$, $2$, $5$, $16$, $61$, $272$, $1385$, $7936$, $50521$
(see~\cite[Sequence A000111(M1492)]{SP}). These numbers are known
as the Euler numbers and have exponential generating function
$\sec x+\tan x$.

A permutation is said to be \emph{up-up} (resp. \emph{up-down},
\emph{down-up}, \emph{down-down}) if it starts with rise (resp.
rise, descent, descent), then descents and rises (resp. descents
and rises, rises and descents, rises and descents) come in turn,
and ends with rise (resp. descent, rise, descent). We denote the
set of all up-up (resp. up-down, down-up, down-down) permutations
of length $n$ by $UU_n$ (resp. $UD_n$, $DU_n$, $DD_n$) for all
$n\geq2$ (for $n=0,1$ we define $UU_n=UD_n=DU_n=DD_n=\vr$).
Clearly, for all $n\geq 2$,
\begin{equation}
A_n=UU_n\cup UD_n. \label{alt1}
\end{equation}

\subsection{Organization of the paper} In this paper we use generating
function techniques to study those alternating (up-up, up-down,
down-up, down-down) permutations that avoid or contain exactly
once $1\mn3\mn2$ and that also avoid or contain exactly once an
arbitrary generalized pattern on $k$ letters.

The paper is organized as follows. The case of alternating (up-up,
up-down, down-up, down-down) permutations that avoid both
$1\mn3\mn2$ and an arbitrary generalized pattern $\tau$ is treated
in Section \ref{sec2}. We derive a simple structure for
alternating (up-up, up-down, down-up, down-down) permutations that
avoid $1\mn3\mn2$. This structure can be used for several
interesting cases, including the classical pattern
$\tau=1\mn2\mn\cdots\mn k$ and $\tau=2\mn3\mn\cdots\mn k\mn1$, the
generalized patterns $\tau=12\mn3\mn\cdots\mn k$ and
$\tau=21\mn3\mn\cdots\mn k$. The case of alternating (up-up,
up-down, down-up, down-down) permutations that avoid $1\mn3\mn2$
and contain exactly once $\tau$ is treated in Section \ref{sec3}.
Here again, we use the structure of alternating (up-up, up-down,
down-up, down-down) permutations that avoid $1\mn3\mn2$ for
several particular cases. The case of alternating (up-up, up-down,
down-up, down-down) permutations that contain $1\mn3\mn2$ exactly
once is treated in Sections \ref{sec4} and \ref{sec5} for avoiding
$\tau$, or containing $\tau$ exactly once, respectively. Finally,
in Section \ref{sec6} we present several directions to extend and
generalize the results of the previous sections, including a
statistics on the set of alternating (up-up, up-down, down-up,
down-down) permutations that avoid $1\mn3\mn2$.

Most of the explicit solutions obtained in Sections
\ref{sec2}--\ref{sec5} involve Chebyshev polynomials of the second
kind.
\section{Avoiding $1\mn3\mn2$ and another pattern $\tau$}\label{sec2}
In this section we consider those alternating (up-up, up-down,
down-up, down-down) permutations in $\SS_n$ that avoid $1\mn3\mn2$
and avoid another generalized pattern $\tau$. We begin by setting
some notation. Let $a_\tau(n)$ denote the number of alternating
permutations in $A_n(1\mn3\mn2,\tau)$, and let
$A_\tau(x)=\sum_{n\gs0}A_\tau(n)x^n$ be the corresponding
generating function.
\begin{center}
\begin{figure}[h]
\epsfxsize=2.0in \epsffile{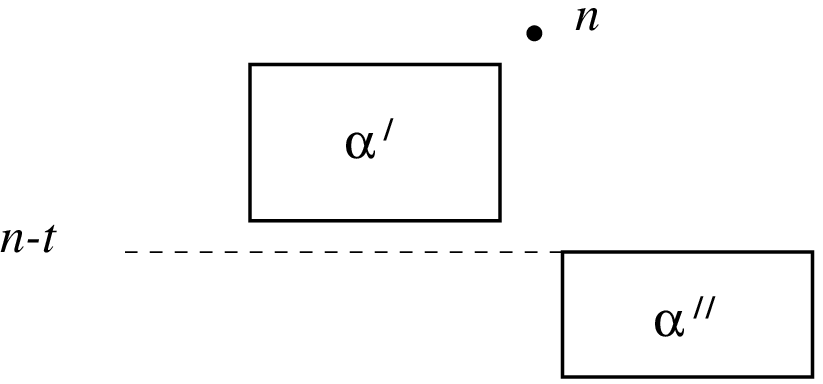} \caption{The block
decomposition for $\alpha\in \SS_n(1\mn3\mn2)$} \label{av132}
\end{figure}
\end{center}
Moreover, we define $ud_\tau(n)$, $uu_\tau(n)$, $du_\tau(n)$, and
$dd_\tau(n)$ as the number of up-up, up-down, down-up, and
down-down permutation in $\SS_n(1\mn3\mn2,\tau)$, respectively. We
denote, the corresponding generating functions by $UD_\tau(x)$,
$UU_\tau(x)$, $DU_\tau(x)$, and $DD_\tau(x)$, respectively. In
this section we describe a method for enumerating alternating
(up-up, up-down, down-up, down-down) permutations that avoid
$1\mn3\mn2$ and another generalized pattern $\tau$ and we use out
this method to enumerate $a_\tau(n)$ for various $\tau$. We begin
with an observation concerning the structure of the alternating
(up-up, up-down, down-up, down-down) permutations that avoid
$1\mn3\mn2$.

Let $\alpha\in \SS_n(1\mn3\mn2)$; by using the block decomposition
approach (see \cite{MV4}) we have one possibility for block
decomposition, as described in Figure~\ref{av132}, if
$\alpha=(\alpha',n,\alpha'')$, then every entry of $\alpha'$ is
greater than every entry of $\alpha''$. Therefore, by using the
block decomposition of $\alpha$, we get the following proposition
which is the base for all the results in this section.

\begin{proposition}\label{av1} Let $\alpha=(\alpha',n,\alpha'')\in
\SS_n(1\mn3\mn2)$.
\begin{enumerate}
\item If $\alpha\in UD_n(1\mn3\mn2)$, then the block
decomposition of $\alpha$ can have one of the following
possibilities:
    \begin{itemize}
    \item[$(1.1)$] $\alpha'$ is an $1\mn3\mn2$-avoiding up-down
    permutation on the letters $n-j+1,n-j+2,\dots,n-1$, and either $\alpha''$
    is an $1\mn3\mn2$-avoiding up-down permutation on the letters
    $1,2,\dots,n-j$ or $\alpha''=1$;

    \item[$(1.2)$] $\alpha'=n-1$ and either $\alpha''$ is an
    $1\mn3\mn2$-avoiding up-down permutation on the letters
    $1,2,\dots,n-j-1$ or $\alpha''=1$.
    \end{itemize}
\medskip
\item If $\alpha\in UU_n(1\mn3\mn2)$, then the block
decomposition of $\alpha$ can have one of the following
possibilities:
    \begin{itemize}
    \item[$(2.1)$] $\alpha'=n-1$ and either $\alpha''=\vr$ or
    $\alpha''$ is an $1\mn3\mn2$-avoiding up-up
    permutation on the letters $1,2,\dots,n-2$;

    \item[$(2.2)$] $\alpha''=\varnothing$ and $\alpha'$ is an $1\mn3\mn2$-avoiding up-down
    permutation on the letters $1,2,\dots,n-1$;

    \item[$(2.3)$] $\alpha'$ is an $1\mn3\mn2$-avoiding up-down permutation on the
    letters $n-j+1,\dots,n-2,n-1$, and $\alpha''$ is an $1\mn3\mn2$-avoiding up-up
    permutation on the letters $1,2,\dots,n-j$.
    \end{itemize}
\medskip
\item If $\alpha\in DU_n(1\mn3\mn2)$, then the block
decomposition of $\alpha$ can have one of the following
possibilities:
    \begin{itemize}
    \item[$(3.1)$] $\alpha'=\vr$, and $\alpha''$ is an $1\mn3\mn2$-avoiding up-up
    permutation on the letters $1,2,\dots,n-1$;

    \item[$(3.2)$] $\alpha'$ is an $1\mn3\mn2$-avoiding down-down permutation on
    the letters $n-j+1,\dots,n-2,n-1$, and either $\alpha''=\vr$
    or $\alpha''$ is an $1\mn3\mn2$-avoiding up-up permutation on the letters
    $1,2,\dots,n-j$.
    \end{itemize}
\medskip
\item If $\alpha\in DD_n(1\mn3\mn2)$, then the block decomposition
of $\alpha$ can have one of the following possibilities:
    \begin{itemize}
    \item[$(4.1)$] $\alpha'=\vr$, and $\alpha''$ is an $1\mn3\mn2$-avoiding up-down
    permutation on the letters $1,2,\dots,n-1$;

    \item[$(4.2)$] $\alpha'$ is an $1\mn3\mn2$-avoiding down-down permutation on the
    letters $n-j+1,\dots,n-2,n-1$, and either $\alpha''=1$ or
    $\alpha''$ is an $1\mn3\mn2$-avoiding up-down permutation on the letters $1,2,\dots,n-j$.
    \end{itemize}
\end{enumerate}
\end{proposition}

\subsection{The pattern $\tau=\vr$} Using Proposition~\ref{av1}, we now
enumerate those alternating (up-up, up-down, down-up, down-down)
permutations in $\SS_n$ that avoid $1\mn3\mn2$.

\begin{theorem}\label{tha0}
The number of alternating permutations of length $n$ that avoid
$1\mn3\mn2$ is given by $C_{[n/2]}$ for any $n\geq0$. Moreover,
for all $n\geq3$,
\begin{enumerate}
\item the number of up-down permutations of length $n$ that avoid $1\mn3\mn2$
is given by $C_{(n-1)/2}$,

\item the number of up-up permutations of length $n$ that avoid $1\mn3\mn2$
is given by $C_{n/2}$,

\item the number of down-up permutation of length $n$ that avoid $1\mn3\mn2$
is given by $C_{(n+1)/2}$,

\item the number of down-down permutations of
length $n$ that avoid $1\mn3\mn2$ is given by $C_{n/2}$,
\end{enumerate}
where $C_m$ is the $m$th Catalan number, such that $C_m$ is
assumed to be $0$ whenever $m$ is a non-integer number.
\end{theorem}
\begin{proof}
An alternating permutation can be either up-down, or up-up, or
down-up, or down-down. If we write an equation for every case in
terms of generating functions and use Proposition~\ref{av1}, then
we get
$$\left\{\begin{array}{l}
UD_\vn(x)=x(UD_\vn(x)+x)^2,\\
DU_\vn(x)=xUU_\vn(x)+xDD_\vn(x)UU_\vn(x)+xDD_\vn(x),\\
DD_\vn(x)=x(x+UD_\vn(x))+xDD_\vn(x)(x+UD_\vn(x)),\\
UU_\vn(x)=x^2(UU_\vn(x)+1)+xUD_\vn(x)+xUD_\vn(x)UU_\vn(x).
\end{array}\right.$$
Therefore, solving the above equations, we get
\begin{equation}
\begin{array}{ll}
UD_\vn(x)=\frac{1-2x^2-\sqrt{1-4x^2}}{2x},& DU_\vn(x)=x\left(\frac{1-\sqrt{1-4x^2}}{2x^2}\right)^2-x,\\
DD_\vn(x)=\frac{1-\sqrt{1-4x^2}}{1+\sqrt{1-4x^2}},& UU_\vn(x)=\frac{1-\sqrt{1-4x^2}}{1+\sqrt{1-4x^2}}.
\end{array}
\label{r0132}\end{equation} Hence, the rest is easy to check by
using $A_\vn(x)=1+x+UD_\vn(x)+UU_\vn(x)$ (see
Equation~\ref{alt1}).
\end{proof}
Using Proposition~\ref{av1}, we now enumerate various sets of
alternating (up-up, up-down, down-up, down-down) permutations in
$\SS_n$ that avoid $1\mn3\mn2$ and nonempty generalized pattern
$\tau$.
\subsection{A classical pattern $\tau=\pk$}
Let us start by the following example.
\begin{example}\label{exaa1}
By definitions we have
$UD_{1}(x)=DD_{1}(x)=UU_{1}(x)=DU_{1}(x)=0$,
$UD_{1\mn2}(x)=UU_{1\mn2}(x)=DU_{1\mn2}(x)=0$, and
$DD_{1\mn2}(x)=x^2$.
\end{example}
The case of varying $k$ is more interesting. As an extension of
Example~\ref{exaa1}, let us consider the case $\tau=\pk$.

\begin{theorem}\label{thaa1}
For all $k\geq 2$,
\begin{enumerate}
\item  $UD_\pk(x)=\dfrac{xU_{k-3}\ttx}{U_{k-1}\ttx}$;

\item $DD_\pk(x)=\dfrac{x^{k-1}+U_{k-3}\ttx}{U_{k-1}\ttx}$;

\item $UU_\pk(x)=\dfrac{U_{k-3}\ttx}{U_{k-1}\ttx}$;

\item $DU_\pk(x)=\dfrac{x^{k-1}+U_{k-3}\ttx}{xU_{k-1}\ttx}-x$.
\end{enumerate}
\end{theorem}
\begin{proof}
Let us verify (1). By using Proposition~\ref{av1}, we have
\begin{equation}
UD_\pk(x)=x(UD_\pk(x)+x)(UD_\pkk(x)+x), \label{ud001}
\end{equation}
equivalently,
$$UD_\pk(x)=\frac{x^2(x+UD_\pkk(x))}{1-x^2-xUD_\pkk(x)}.$$
Now, assume that (1) holds for $k-1$. So, by Identity~\ref{rcheb}
we have that $$UD_\pkk(x)=x^3R_{k-3}(x^2)R_{k-2}(x^2).$$
Therefore,
$$UD_\pk(x)=\frac{x^2(x+x^3R_{k-3}(x^2)R_{k-2}(x^2))}{1-x^2-x^4R_{k-3}(x^2)R_{k-2}(x^2)},$$
and using Identity~\ref{rk}, we get
$$UD_\pk(x)=x^3R_{k-2}(x^2)R_{k-1}(x^2).$$ So, by
Equation~\ref{rcheb}, (1) holds for $k$. Hence, by induction on
$k$ together with Example~\ref{exaa1} we get the desired result as
claimed in (1).

Let us verify (2). By using Proposition~\ref{av1}, we get
$$DD_\pk(x)=x(x+UD_\pk(x))(1+DD_\pkk(x)).$$
Using (1), Equation~\ref{rk}, and induction on $k$ together with
Example~\ref{exaa1} we get the desired result as claimed in (2).

Similarly, we have
$$\begin{array}{l}
UU_\pk(x)=x(x+UD_\pkk(x))(1+UU_\pk(x)),\\
DU_\pk(x)=xUU_\pk(x)+xDD_\pkk(x)(UU_\pk(x)+1),
\end{array}$$
hence, by using Equation~\ref{rk} and the cases (1)--(2), we get
the desired result as claimed in (3) and (4).
\end{proof}

As a corollary to Theorem~\ref{thaa1} and Equation~\ref{alt1}, we
get the number of alternating permutations in
$A_n(1\mn3\mn2,\pk)$.

\begin{corollary}\label{ccin1}
For all $k\geq 2$,
$$A_\pk(x)=\dfrac{(1+x)U_{k-2}\ttx}{xU_{k-1}\ttx}=(1+x)R_{k-1}(x^2).$$
\end{corollary}

Corollary~\ref{ccin1}, for $k=5$, yields that the number of
alternating permutations in $A_n(1\mn3\mn2,1\mn2\mn3\mn4\mn5)$ is
given by $F_{[(n+2)/2]}$, where $F_m$ is the $m$th Fibonacci
number.
\subsection{A classical pattern $\tau=2\mn3\mn\cdots\mnk\mn1$}
By Proposition \ref{av1}, we have that $UD_{2\mn3\mn1}(x)=0$ and
$DD_{2\mn3\mn1}(x)=x^2$. The case of varying $k$ is more
interesting. As an extension of this example, let us consider the
case $\tau=2\mn3\mn\cdots\mnk\mn1$.

\begin{theorem}\label{thab1}
For all $k\geq 3$,
\begin{enumerate}
\item
$UD_{2\mn3\mn\cdots\mnk\mn1}(x)=\dfrac{xU_{k-4}\ttx}{U_{k-2}\ttx}$,

\item
$DD_{2\mn3\mn\cdots\mnk\mn1}(x)=\dfrac{x^{k-2}+U_{k-4}\ttx}{U_{k-2}\ttx}$,

\item
$UU_{2\mn3\mn\cdots\mnk\mn1}(x)=\dfrac{U_{k-3}^2\ttx}{U_{k-2}^2\ttx}$,

\item
$DU_{2\mn3\mn\cdots\mnk\mn1}(x)=\dfrac{x}{U_{k-2}^2\ttx}\biggl(x^{k-3}(U_{k-3}\ttx+xU_{k-2}\ttx)+2U_{k-3}^2\ttx+U_{k-4}^2\ttx-2\biggr)$.
\end{enumerate}
\end{theorem}
\begin{proof}
Proposition \ref{av1}(1) yields
$$UD_{2\mn3\mn\cdots\mnk\mn1}(x)=x(UD_{2\mn3\mn\cdots\mnk\mn1}(x)+x)(UD_{1\mn2\mn\cdots\mn(k-2)}(x)+x).$$
Using Equation~\ref{rk} and Recurrence~\ref{reccheb}, together
with Theorem~\ref{thaa1}(1), we have that case (1) holds.

Proposition \ref{av1}(4) yields
$$DD_{2\mn3\mn\cdots\mnk\mn1}(x)=x(x+UD_{2\mn3\mn\cdots\mnk\mn1}(x))(1+UD_{1\mn2\mn\cdots\mn(k-2)}x)).$$
By using the first case (1) and Theorem~\ref{thaa1}(1), together
with Equation~\ref{rk} and Recurrence~\ref{reccheb}, we get that
case (2) holds.

Proposition \ref{av1}(2) yields
$$UU_{2\mn3\mn\cdots\mnk\mn1}(x)=x^2(1+UU_{2\mn3\mn\cdots\mnk\mn1}(x))+xUD_{2\mn3\mn\cdots\mnk\mn1}(x)+xUD_{1\mn2\mn\cdots\mn(k-2)}(x)UU_{2\mn3\mn\cdots\mnk\mn1}(x)).$$
By using cases (1) and (2), together with Theorem~\ref{thaa1},
Equation~\ref{rk}, and Recurrence~\ref{reccheb}, we have that case
(3) holds.

Proposition \ref{av1}(3) yields
$$DU_{2\mn3\mn\cdots\mnk\mn1}(x)=xUU_{2\mn3\mn\cdots\mnk\mn1}(x)+xDD_{1\mn2\mn\cdots\mn(k-2)}(x)UU_{2\mn3\mn\cdots\mnk\mn1}(x)+xDD_{2\mn3\mn\cdots\mnk\mn1}(x).$$
By using cases (3) and (4), Theorem~\ref{thaa1}(2) and
Recurrence~\ref{reccheb}, we get that case (4) holds.
\end{proof}

Theorem~\ref{thab1} yields an explicit formula to the generating
function for the number of alternating permutations in
$A_n(1\mn3\mn2,2\mn3\mn\cdots\mnk\mn1)$.

\begin{corollary}\label{cthab}
For all $k\geq 3$,
$$A_{2\mn3\mn\cdots\mnk\mn1}(x)=\frac{(1+x)U_{k-3}^2\ttx-x}{U_{k-2}^2\ttx}.$$
\end{corollary}
For example, by using Corollary~\ref{cthab} we get that for all
$n\geq1$,
$$A_{2\mn3\mn4\mn5\mn6\mn1}(2n)=\frac{7}{10}nL_{2n}-\frac{1}{10}(15n-4)F_{2n}\mbox{
and }A_{2\mn3\mn4\mn5\mn6\mn1}(2n+1)=F_{2n-1},$$ where $F_m$ is
the $m$th Fibonacci number and $L_m$ is the $m$th Lucas number.
\subsection{A generalized patterns $\tau\mn3\mn\cdots\mn k$} In the
current subsection we interesting in alternating permutations with
two restrictions, the first one is $1\mn3\mn2$ (classical pattern
$132$) and the second one is $\tau\mn3\mn\cdots\mn k$ (generalized
pattern), where $\tau=12$, $\tau=21$, $\tau=1\mn2$, or
$\tau=2\mn1$. Using the same arguments as in the proof of
Theorem~\ref{thaa1}, we get

\begin{theorem}\label{thac1}
Let $\tau\in\{12, 21, 1\mn2, 2\mn1\}$. For all $k\geq 2$,
\begin{enumerate}
\item
$UD_{\tau\mn3\mn\cdots\mnk}(x)=\dfrac{xU_{k-3}\ttx}{U_{k-1}\ttx}$,

\item
$DD_{\tau\mn3\mn\cdots\mnk}(x)=\dfrac{x^{k-1}+U_{k-3}\ttx}{U_{k-1}\ttx}$,

\item
$UU_{\tau\mn3\mn\cdots\mnk}(x)=\dfrac{U_{k-3}\ttx}{U_{k-1}\ttx}$,

\item $DU_{\tau\mn3\mn\cdots\mnk}(x)=\dfrac{x^{k-1}+U_{k-3}\ttx}{xU_{k-1}\ttx}-x$.
\end{enumerate}
\end{theorem}
The above theorem together with Equation~\ref{alt1} yield that the
number of alternating permutations in $A_n(1\mn3\mn2,
\tau\mn3\mn\cdots\mnk)$, where $\tau\in\{12,21,1\mn2,2\mn1\}$.

\begin{corollary}\label{ccaa1}
Let $\tau\in\{12, 21, 1\mn2, 2\mn1\}$. For any $k\geq 2$,
$$A_{\tau\mn3\mn\cdots\mnk}(x)=\dfrac{(1+x)U_{k-2}\ttx}{xU_{k-1}\ttx}=(1+x)R_{k-1}(x^2).$$
\end{corollary}

Corollary~\ref{ccaa1} suggests that there should exist a bijection
between the sets $A_n(1\mn3\mn2,\pk)$ and
$A_n(1\mn3\mn2,2\mn1\mn3\mn\cdots\mnk)$. However, we failed to
produce such a bijection, and finding it remains a challenging
open question.
\section{Avoiding $1\mn3\mn2$ and containing another
pattern}\label{sec3} In this section we consider those alternating
(up-up, up-down, down-up, down-down) permutations in $\SS_n$ that
avoid $1\mn3\mn2$ and contain another generalized pattern $\tau$
exactly once. We begin by setting some notation. Let
$a_{\tau;r}(n)$ be the number of alternating permutations in
$A_n(1\mn3\mn2)$ that contain $\tau$ exactly $r$ times. Moreover,
we denote the number of up-up (resp. up-down, down-up, down-down)
permutations in $\SS_n(1\mn3\mn2)$ that contain $\tau$ exactly $r$
times by $uu_{\tau;r}(n)$ (resp. $ud_{\tau;r}(n)$,
$du_{\tau;r}(n)$, $dd_{\tau;r}(n)$). We denote the corresponding
generating function by $A_{\tau;r}(x)$ (resp. $UU_{\tau;r}(x)$,
$UD_{\tau;r}(x)$, $DU_{\tau;r}(x)$, $DD_{\tau;r}(x)$). In this
section we use out Proposition~\ref{av1} and generating function
techniques to enumerate $a_{\tau;r}(n)$, $ud_{\tau;r}(n)$,
$uu_{\tau;r}(x)$, $du_{\tau;r}(n)$, and $dd_{\tau;r}(n)$ for
various cases of $\tau$.
\subsection{A classical pattern $\tau=\pk$}
\begin{theorem}\label{thc1}
For all $k\geq 2$,
\begin{enumerate}
\item $UD_{\pk;1}(x)=\dfrac{x}{U_{k-1}^2\ttx}$;

\item $DD_{\pk;1}(x)=\dfrac{1}{U_{k-1}\ttx}\sum\limits_{m=0}^{k-2}\dfrac{x^{k-1-m}(x^{m+1}+U_{m-1}\ttx)}{U_m\ttx
U_{m+1}\ttx}$;

\item $UU_{\pk;1}(x)=\dfrac{1}{U_{k-1}^2\ttx}$;

\item $DU_{\pk;1}(x)=\dfrac{1}{xU_{k-1}\ttx}\sum\limits_{m=0}^{k-2}\dfrac{x^{k-1-m}(x^{m+1}+U_{m-1}\ttx)}{U_m\ttx
U_{m+1}\ttx}$.
\end{enumerate}
\end{theorem}
\begin{proof}
To verify (1), as consequence of Proposition~\ref{av1}, we have
three possibilities for the block decomposition for an arbitrary
up-down permutation $\pi=(\alpha,n,\beta)\in UD_n(1\mn3\mn2)$ that
contain $\pk$ exactly once. The first possibility is
$\alpha=(n-1)$ and $\beta\in UD_{n-2}(1\mn3\mn2)$ contains $\pk$
exactly once. The second possibility is $\beta=(1)$ and $\alpha$
is up-down permutation of the numbers $2,3,\cdots,n-1$ which
contains $\pkk$ exactly once. The third possibility satisfies the
condition that every entry of $\alpha$ is greater than every entry
of $\beta$, where $\alpha$ and $\beta$ are up-down permutations on
the numbers $n-j+1,n-j+2,\cdots,n-1$ and $1,2,\dots,n-j$
(respectively) such that either $\alpha$ contains $\pkk$ exactly
once and $\beta$ avoids $\pk$, or $\alpha$ avoids $\pkk$ and
$\beta$ contains $\pk$ exactly once. So, the first, the second,
and the third cases above give contribution $x^2UD_{\pk;1}(x)$,
$x^2UD_{\pkk;1}(x)$,
$xUD_{\pkk;1}(x)UD_{\pk;0}(x)+xUD_{\pkk;0}(x)UD_{\pk;1}(x)$,
respectively. Therefore,
$$\begin{array}{l}
UD_{\pk;1}(x)=x^2UD_{\pk;1}(x)+x^2UD_{\pkk;1}(x)+\\
\qquad\qquad\qquad+xUD_{\pkk;1}(x)UD_{\pk;0}(x)+xUD_{\pkk;0}(x)UD_{\pk;1}(x).
\end{array}$$
Using Theorem~\ref{thaa1} together with Equation~\ref{rk} and
Recurrence~\ref{reccheb}, we get
$$UD_{\pk;1}(x)=xR_{k-1}^2(x^2)UD_{\pkk;1}(x).$$
Besides, by definitions we have that $UD_{1\mn2;1}(x)=x^2$. Hence,
by Identity~\ref{rcheb} and by induction,  case (1) holds.

Let us verify (2). Using the same arguments as in the proof of
(1), we get
$$\begin{array}{l}
DD_{\pk;1}(x)=xUD_{\pk;1}(x)+xDD_{\pkk;1}(x)UD_{\pk;0}(x)+\\
\qquad\qquad\qquad+xDD_{\pkk;0}(x)UD_{\pk;1}(x)+x^2DD_{\pkk;1}(x),
\end{array}$$ together with $DD_{1\mn2;1}(x)=x^4$ (by the definitions).
Hence, by using (1), Theorem~\ref{thaa1} and induction on $k$, we
get the desired result as claimed in (2).

Let us verify (3). Using the same arguments as in the proof of
(1), we get
$$\begin{array}{l}
UU_{\pk;1}(x)=x^2UU_{\pk;1}(x)+xUD_{\pkk;1}(x)+\\
\qquad\qquad\qquad+xUD_{\pkk;0}(x)UU_{\pk;1}(x)+xUD_{\pkk;1}(x)UU_{\pk;0}(x),
\end{array}$$ equivalently,
$$UU_{\pk;1}(x)=\frac{xUU_{\pk;0}(x)UD_{\pkk;1}(x)}{1-x^2-xUD_{\pkk;0}(x)}.$$
Hence, by using (1), Equation~\ref{rk} and
Recurrence~\ref{reccheb}, case (3) holds.

Similarly to the cases above, case (4) holds.
\end{proof}

As a corollary to Theorem~\ref{thc1} and Equation~\ref{alt1}, we
get the number of alternating permutations in $\SS_n(1\mn3\mn2)$
that contain $\pk$ exactly once.

\begin{corollary}\label{ddaa1}
For all $k\geq 3$,
$$A_{\pk;1}(x)=\dfrac{1+x}{U_{k-1}^2\ttx}.$$
\end{corollary}

For example, the number of alternating permutations in
$A_n(1\mn3\mn2)$ that contain $1\mn2\mn3\mn4\mn5$ exactly once is
given by $\sum_{j=1}^{[n/2]-3} F_{2j}F_{2[n/2]-4-2j}$ for all
$n\geq6$, where $F_m$ is the $m$th Fibonacci number.
\subsection{A generalized patterns $\tau\mn3\mn\cdots\mn k$} In the
current subsection we interesting in alternating permutations with
two restrictions, the first one is $1\mn3\mn2$ (classical pattern
$132$) and the second one is $\tau\mn3\mn\cdots\mn k$ (generalized
pattern), where $\tau=12$, $\tau=21$, $\tau=1\mn2$, or
$\tau=2\mn1$. Using the same arguments as in the proof of
Theorem~\ref{thc1}, we get

\begin{theorem}\label{thc2}
Let $\tau\in\{12, 21, 1\mn2\}$. For all $k\geq 3$,
\begin{enumerate}
\item $UD_{\tau\mn3\mn\cdots\mnk;1}(x)=\frac{x}{U_{k-1}^2\ttx}$;

\item $UU_{\tau\mn3\mn\cdots\mnk;1}(x)=\frac{1}{U_{k-1}^2\ttx}$,

\item $UD_{2\mn1\mn3\cdots\mnk;1}(x)=UU_{2\mn1\mn3\cdots\mnk;1}(x)=0$.
\end{enumerate}
\end{theorem}

As a corollary to Theorem~\ref{thc2} and Equation~\ref{alt1} we
get the number of alternating permutations in $A_n(1\mn3\mn2)$
that contain $\tau\mn3\mn\cdots\mnk$ exactly once, where
$\tau\in\{12,21,1\mn2,2\mn1\}$.

\begin{corollary}\label{ccaa2}
Let $\tau\in\{12, 21, 1\mn2\}$. For $k\geq 3$,
$$A_{\tau\mn3\mn\cdots\mnk;1}(x)=\dfrac{(1+x)U_{k-2}\ttx}{xU_{k-1}\ttx}=(1+x)R_{k-1}(x^2).$$
Besides, $A_{2\mn1\mn3\mn\cdots\mnk;1}(x)=0$ for all $k\geq 3$.
\end{corollary}

Corollary~\ref{ccaa2} suggests that there should exist a bijection
between the sets of alternating permutations of length $n$ that
avoid $1\mn3\mn2$ and contain $\pk$ exactly once and the set of
alternating permutations of length $n$ that avoid $1\mn3\mn2$ and
contain $2\mn1\mn3\mn\cdots\mnk$ exactly once. However, we failed
to produce such a bijection, and finding it remains a challenging
open question.

Again, by the same arguments as in the proof of
Theorem~\ref{thc1}, with using Proposition~\ref{av1}, we get the
following result.

\begin{theorem}
For all $k\geq 3$,
\begin{enumerate}
\item $DD_{\pk;1}(x)=\frac{x}{U_{k-1}\ttx}\sum_{j=0}^{k-2}
\frac{x^{k-2-j}(x^{j+1}+U_{j-1}\ttx)}{U_j\ttx U_{j+1}\ttx}$,

\item $DD_{12\mn3\mn4\mn\cdots\mnk;1}(x)=\frac{x}{U_{k-1}\ttx}\left(
x^k+\sum_{j=0}^{k-2} \frac{x^{k-2-j}(x^{j+1}+U_{j-1}\ttx)}{U_j\ttx
U_{j+1}\ttx} \right)$,

\item
$DD_{2\mn1\mn3\mn4\mn\cdots\mnk;1}(x)=\frac{x^{k-1}}{U_{k-1}\ttx}$,

\item $DD_{21\mn3\mn4\mn\cdots\mnk;1}(x)=\frac{1}{U_{k-1}^2\ttx}$.
\end{enumerate}
Moreover, for all $k\geq3$,
$$DU_{\tau\mn3\mn4\mn\cdots\mnk;1}(x)=\frac{1}{x}DD_{\tau\mn3\mn4\mn\cdots\mnk;1}(x),$$
where $\tau\in\{12,21,1\mn2,2\mn1\}$.
\end{theorem}
\subsection{A three letter generalized pattern without dashes} In this
subsection we find an explicit expression for $UD_{\tau;r}(x)$,
$UU_{\tau;r}(x)$, $DU_{\tau;r}(x)$, $DD_{\tau;r}(x)$, and
$A_{\tau;r}(x)$, where $\tau=123$, $213$, $231$, $312$, $321$ are
generalized pattern without dashes.

\begin{theorem}\label{thud3}
For all $r\geq0$,
\begin{enumerate}
\item $UD_{123;r}(x)=UD_{321}(x;r)=0$,

\item $UD_{213;r}(x)=\sum_{n\geq r+1}\frac{(r+1)}{n(n-r)}\binom{n}{r+1}^2
x^{2n+1}$,

\item $UD_{231;r}(x)=\frac1{r+1}\binom{2r}{r}x^{2r+1}$ for $r\geq1$,

\item $UD_{312;r}(x)=\sum_{n\geq r+1}\frac{(r+1)}{n(n-r)}\binom{n}{r+1}^2
x^{2n+1}$.
\end{enumerate}
\end{theorem}
\begin{proof}
First of all, let us define $UD_{\tau}(x,q)=\sum_{r\geq0}
UD_{\tau;r}(x)q^r$ (similarly, we define $UU_{\tau}(x,q)$,
$DU_{\tau}(x,q)$, and $DD_{\tau}(x,q)$).

The first case (1) holds by the definitions. To verify (2), let us
use the same arguments as in the proof of Theorem~\ref{thaa1}. So,
we have
$$UD_{213}(x,q)=x^2(x+UD_{213}(x,q))+xqUD_{213}(x,q)(x+UD_{213}(x,q)),$$
equivalently,
$$UD_{213}(x,q)=\frac{1-x^2(1+q)-\sqrt{x^4(1-q)^2-2x^2(q+1)+1}}{2xq}.$$
The rest is easy to check.

Similarly, we have that
$$UD_{231}(x,q)=x^2q(x+UD_{231}(x,q))+xqUD_{231}(x,q)(x+UD_{231}(x)),$$
equivalently,
$$UD_{231}(x,q)=\frac{x^3q}{1-2x^2q-xqUD_{231}(x,q)}.$$
The rest is easy to check.

Finally, by the reversal symmetric operation (defined by
$(\pi_1,\dots,\pi_n)\mapsto(\pi_n,\dots,\pi_1)$) we get that
$UD_{213}(x;r)=UD_{312}(x;r)$ for all $r\geq0$.
\end{proof}

Similarly to Theorem~\ref{thud3} we get the following result.

\begin{theorem}\label{thuu3}
Let $r\geq0$.
\begin{enumerate}
\item  $UU_{123}(x;r)=UU_{321}(x;r)=0$;

\item $UU_{213}(x;r)=\sum_{n\geq r+1}\frac{(r+1)}{n(n-r)}\binom{n}{r+1}^2 x^{2n}$;

\item  $UD_{231}(x;r)=\frac1{r+1}\binom{2r}{r}x^{2r+2}$ for $r\geq1$;

\item  $UU_{312}(x;r)=\sum_{n\geq r+1}\frac{(r+1)}{n(n-r)}\binom{n}{r+1}^2 x^{2n}$.
\end{enumerate}
\end{theorem}

As a corollary to Theorem~\ref{thud3} and Theorem~\ref{thuu3}, we
have the following result.

\begin{theorem}\label{thalt3}
Let $r\geq0$.
\begin{enumerate}
\item $A_{123}(x;r)=A_{321}(x;r)=0$;

\item $A_{213}(x;r)=A_{312}(x;r)=\sum_{n\geq r+1}\frac{(r+1)}{n(n-r)}\binom{n}{r+1}^2 (x^{2n}+x^{2n+1})$;

\item $A_{231}(x;r)=C_r(x^{2r+2}+x^{2r+1})$ for $r\geq1$.
\end{enumerate}
\end{theorem}

\begin{center}
\begin{figure}[h]
\epsfxsize=2.2in \epsffile{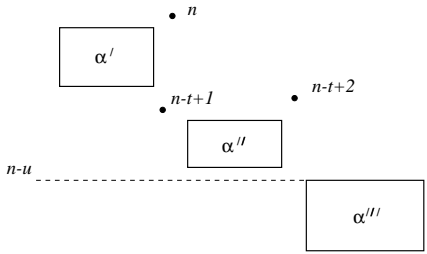} \caption{The second
contribution for the block decomposition of a permutation that
contain $1\mn3\mn2$ exactly once}\label{co132}
\end{figure}
\end{center}
\section{Containing $1\mn3\mn2$ exactly once and avoiding arbitrary
pattern}\label{sec4} In this section we consider those alternating
(up-up, up-down, down-up, down-down) permutations in $\SS_n$ that
contain $1\mn3\mn2$ exactly once and avoid an arbitrary
generalized pattern $\tau$. We begin by setting some notation. Let
$A^1_{\tau}(x)$ be the generating function for the number of
alternating permutations in $A_n(\tau)$ that contain $1\mn3\mn2$
exactly once. More precisely, We denote the generating function
for number of up-up (resp. up-down, down-up, down-down)
permutations in $\SS_n(\tau)$ that contain $1\mn3\mn2$ exactly
once by $UD^1_{\tau}(x)$ (resp. $UU^1_{\tau}(x)$,
$DU^1_{\tau}(x)$, $DD^1_{\tau}(x)$). In this section we describe a
method for enumerating alternating (up-up, up-down, down-up,
down-down) permutations that contain $1\mn3\mn2$ exactly once and
avoid another generalized pattern. We begin with observation
concerning the structure of the permutations that contain
$1\mn3\mn2$ exactly once.

Using the block decomposition approach (see~\cite{MV4}) we get two
possibilities for the block decomposition for permutation $\pi$
containing $1\mn3\mn2$ exactly once. The first one is
$\alpha=(\alpha',n,\alpha'')$, where every entry of $\alpha'$ is
greater than every entry of $\alpha''$, and the second one is
$\alpha=(\alpha',n-t+1,n,\alpha'',n-t+2,\alpha''')$, where every
entry of $\alpha'$ is greater than $n-t+2$, $n-t+1$ is greater
than every entry of $\alpha''$, and every entry of $\alpha''$ is
greater than every entry of $\alpha'''$, as described in
Figures~\ref{av132} and \ref{co132}. Therefore, by using the block
decomposition of $\alpha$, we get the following proposition which
is the base for all the results in this section.

\begin{proposition}\label{cv1}
Let $\alpha\in\SS_n$ contains $1\mn3\mn2$ exactly once. Then the
block decomposition of $\alpha$ can have one of the following
three forms:

(i) there exists $t$ such that $\alpha=(\alpha',n,\alpha'')$,
where $\alpha'$ is a permutation of $n-1,n-2,\dots,n-t+1$ which
contains $1\mn3\mn2$ exactly once, and $\alpha''$ is a permutation
of $1,2,\dots,n-t$ which avoids $1\mn3\mn2$;

(ii) there exists $t$ such that $\alpha=(\alpha',n,\alpha'')$,
where $\alpha'$ is a permutation of $n-1,n-2,\dots,n-t+1$ which
avoids $1\mn3\mn2$, and $\alpha''$ is a permutation of
$1,2,\dots,n-t$ which contains $1\mn3\mn2$ exactly once;

(iii) there exist $t,u$ such that
$\alpha=(\alpha',n-t+1,n,\alpha'',n-t+2,\alpha''')$, where
$\alpha'$ is a permutation of $n-1,n-2,\dots,n-t+3$ which avoids
$1\mn3\mn2$, $\alpha''$ is a permutation of
$n-t,n-t-1,\dots,n-u+1$ which avoids $1\mn3\mn2$, and $\alpha'''$
is a permutation of $1,2,\dots,n-u$ which avoids $1\mn3\mn2$.
\end{proposition}

As a remark, using the proposition above for an arbitrary
permutation $\alpha$ which contains $1\mn3\mn2$ exactly once
together with the definitions of up-up, up-down, down-up,
down-down permutations we can present all the possibilities of the
block decomposition for any arbitrary up-up, up-down, down-up,
down-down permutation, respectively.

\subsection{The pattern $\tau=\vr$} Using Proposition~\ref{cv1},
we now enumerate those alternating (up-up, up-down, down-up,
down-down) permutations in $\SS_n$ that contain $1\mn3\mn2$
exactly once.

\begin{theorem}\label{thba1} We have
\begin{enumerate}
\item The generating function for the number of up-down
permutations in $\SS_n$ containing $1\mn3\mn2$ exactly once is
given by
$$UD^1_\vr(x)=\frac{x(1-\sqrt{1-4x^2})}{1-4x^2+\sqrt{1-4x^2}}.$$
In other words, the number of up-down permutations in $\SS_n$
containing $1\mn3\mn2$ exactly once is given by
$\binom{n-1}{(n-3)/2}$,

\item The generating function for the number of up-up
permutations in $\SS_n$ containing $1\mn3\mn2$ exactly once is
given by
$$UU^1_\vr(x)=\frac{x^2-1}{x^2}+\frac{1-3x^2}{x^2\sqrt{1-4x^2}}.$$
In other words, the number of up-up permutations in $\SS_n$
containing $1\mn3\mn2$ exactly once is given by
$2\binom{n-1}{(n-4)/2}$,

\item The generating function for the number of down-down
permutations in $\SS_n$ containing $1\mn3\mn2$ exactly once is
given by
$$DD^1_\vr(x)=\frac{x^2-1}{x^2}+\frac{1-3x^2}{x^2\sqrt{1-4x^2}}.$$
In other words, the number of down-down permutations in $\SS_n$
containing $1\mn3\mn2$ exactly once is given by
$2\binom{n-1}{(n-4)/2}$,

\item The generating function for the number of down-up
permutations in $\SS_n$ containing $1\mn3\mn2$ exactly once is
given by
$$DU^1_\vr(x)=\frac{2x^4+4x^2-3}{2x^3}+\frac{3-10x^2}{2x^3\sqrt{1-4x^2}}.$$
In other words, the number of down-up permutations in $\SS_n$
containing $1\mn3\mn2$ exactly once is given by
$\frac{3}{2}\binom{n+3}{(n+3)/2}-5\binom{n+1}{(n+1)/2}$,
\end{enumerate}
where $\binom ab$ is assumed to be $0$ whenever $a<b$, $b<0$, or
$b$ is a non-integer number.
\end{theorem}
\begin{proof}
To verify (1), let us consider the three possibilities of the
block decompositions, as described in Proposition~\ref{cv1}, for
an arbitrary permutation $\pi\in \SS_n$ that contain $1\mn3\mn2$
exactly. The first (or the second) contribution of the block
decomposition above for up-down permutations gives
            $$xUD^1_\vr(x)(x+UD_\vr(x)),$$
and the third contribution of the block decomposition above for
up-down permutations gives
        $$x^3\bigl(UU_\vr(x)(x+UD_\vr(x))^2+(x+UD_\vr(x))^2+UU_\vr(x)+1\bigr).$$
Therefore, by using the above three contributions, we have
$$\begin{array}{l}
UD^1_\vr(x)=2xUD^1_\vr(x)(x+UD_\vr(x))+x^3\bigl(UU_\vr(x)(x+UD_\vr(x))^2+(x+UD(x))^2+UU_\vr(x)+1\bigr).
\end{array}$$ Hence, by solving the above equation together with
Theorem~\ref{tha0}, we get the claimed in (1).

Similarly, let us verify (2). The first contribution of the block
decomposition for up-up permutations gives
$xUD^1_\vr(x)(UU_\vr(x)+1)$, the second contribution of the block
decomposition for up-up permutations gives
$x(UD_\vr(x)+x)UU^1_\vr(x)$, and the third contribution of the
block decomposition for up-up permutations gives
$x^3(UU_\vr(x)+1)^2(UD_\vr(x)+x)$. Therefore, by using the three
three contributions above, we get
$$UU^1_\vr(x)=xUD^1_\vr(x)(UU_\vr(x)+1)+x(UD_\vr(x)+x)UU^1_\vr(x)+x^3(UU_\vr(x)+1)^2(UD_\vr(x)+x).$$
Hence, by using (1) together with Theorem~\ref{tha0} we get the
claimed in (2).

Similarly to (1) and (2), we can verify cases (3) and (4).
\end{proof}

As a corollary to Theorem~\ref{thba1} together with
Identity~\ref{alt1} we have

\begin{corollary}
The number of alternating permutations in $A_n$ containing
$1\mn3\mn2$ exactly once is given by
$$\binom{n-1}{(n-3)/2}+\binom{n-1}{(n-4)/2},$$
where $\binom ab$ is assumed to be $0$ whenever $a<b$, $b<0$, or
$b$ is a non-integer number.
\end{corollary}

Using Proposition~\ref{cv1} for alternating permutations, we now
enumerate various sets of alternating (up-up, up-down, down-up,
down-down) permutations in $\SS_n$ that contain $1\mn3\mn2$
exactly once and avoid nonempty generalized pattern $\tau$.

\subsection{A classical pattern $\tau=\pk$}

\begin{theorem}\label{thbb1} For all $k\geq 3$,
\begin{enumerate}
\item $UD^1_\pk(x)=\dfrac{1}{U_{k-1}^2\ttx}\sum\limits_{j=0}^{k-3}U_{j+1}\ttx U_j\ttx$;

\item $UU^1_\pk(x)=\dfrac{1}{xU_{k-1}^2\ttx}\left(xU_{k-3}^2\ttx +\sum\limits_{j=0}^{k-4} U_j\ttx U_{j+1}\ttx\right)$;

\item The generating function $DD^1_\pk(x)$ is given by
$$\frac{1}{U_{k-1}\ttx}\sum\limits_{i=0}^{k-2}x^i\left[
\frac{x^{k-1-i}+U_{k-3-i}\ttx}{U_{k-1-i}\ttx
U_{k-2-i}\ttx}\left(\sum\limits_{j=0}^{k-2-i}U_j\ttx
U_{j+1}\ttx\right)-xU_{k-2-i}\ttx \right];$$

\item The generating function $DU^1_\pk(x)$ is given by
$$\frac{x^{k-1}+U_{k-3}\ttx}{U_{k-2}\ttx}UU_\pk^1(x)+\frac{U_{k-2}\ttx}{U_{k-1}\ttx}DD_{k-1}^1(x)+
\frac{xU_{k-3}\ttx(x^{k-2}+U_{k-4}\ttx)}{U_{k-2}\ttx
U_{k-1}\ttx}.$$
\end{enumerate}
\end{theorem}
\begin{proof}
Using the same arguments as in the proof of Theorem~\ref{thba1},
we get the following. The first contribution of the block
decomposition for up-down permutations gives
$$x^2UD^1_\pk(x)+xUD^1_\pkk(x)(x+UD_\pk(x)),$$
the second contribution of the block decomposition for up-down
permutations gives  $$xUD_\pkk(x)UD^1_\pk(x),$$ and the third
contribution of the block decomposition for up-down permutations
gives
$$x^3(1+UU_\pkk(x))\bigl( (x+UD_\pkk(x))(x+UD_\pk(x))+1\bigr).$$
Therefore, by using the three contributions above together with
Theorem~\ref{thaa1} and Recurrence~\ref{reccheb}, we have
$$UD^1_\pk(x)=\frac{U_{k-2}\ttx U_{k-3}\ttx}{U_{k-1}^2\ttx}
+\frac{U_{k-2}^2\ttx}{U_{k-3}^2\ttx}UD^1_\pkk(x).$$ Hence, by
induction on $k$ together with $UD_{1\mn2}^1(x)=0$ (by the
definitions) we get the desired result as claimed in (1).

Similarly, let us verify (2). The first contribution of the block
decomposition for up-up permutations gives
$$x^2UU^1_\pk(x)+xUD^1_\pkk(x)(1+UU_\pk(x)),$$
the second contribution of the block decomposition for up-up
permutations gives $$xUD_\pkk(x)UU^1_\pk(x),$$ and the third
contribution of the block decomposition for up-up permutations
gives
$$x^3(x+UD_\pkk(x))(1+UU_\pkk(x))(1+UU_\pk(x))..$$
Therefore, by using the three contributions above together with
case (1), Theorem~\ref{thaa1} and Recurrence~\ref{reccheb}, we get
the desired result as claimed in (2).

Similarly, we have
$$\begin{array}{l}
DD_\pk^1(x)=x(1+DD_\pkk(x))UD_\pk^1(x)+x(x+UD_\pk(x))DD_\pkk^1(x)\\
\qquad\qquad\qquad+x^3(x+UD_\pkk(x))(1+(x+UD_\pkk(x))(x+UD_\pk(x))),\\
\\
DU_\pk^1(x)=x(1+DD_{\pkk}(x))UU_\pk^1(x)+x(1+UU_{\pk}(x))DD_{\pkk}^1(x)\\
\qquad\qquad\qquad+x^3(x+DU_\pkk(x))(x+UD_\pkk(x))(1+UU_\pk(x)).
\end{array}$$
Using Theorem~\ref{thaa1}, cases (1) and (2), and
Recurrence~\ref{reccheb}, we get (3).
\end{proof}

As a corollary to Theorem~\ref{thbb1} together with using
Identity~\ref{alt1}, we have

\begin{corollary}
For all $k\geq3$, the generating function $A^1_\pk(x)$ given by
$$\frac{1}{xU_{k-1}^2\ttx}\left[
xU_{k-3}\ttx\left(U_{k-3}\ttx+U_{k-4}\ttx\right)
+(1+x)\sum\limits_{m=0}^{k-4}U_{m+1}\ttx U_m\ttx\right].$$
\end{corollary}

\subsection{A generalized patterns $\tau=12\mn3\mn\cdots\mnk$ or $\tau=21\mn3\mn\cdots\mnk$}
Using the same arguments as in the proof of Theorem~\ref{thbb1},
we get

\begin{theorem}\label{thbb3} For all $k\geq 3$,
$$\begin{array}{lll}
&UD^1_{12\mn3\mn\cdots\mnk}(x)=UD^1_\pk(x),\qquad&UU^1_{12\mn3\mn\cdots\mnk}(x)=UU^1_\pk(x)\\
&DU^1_{12\mn3\mn\cdots\mnk}(x)=DU^1_\pk(x),
&DD^1_{12\mn3\mn\cdots\mnk}(x)=DD^1_\pk(x).
\end{array}$$
Moreover, for all $k\geq3$,
$$A^1_{\tau\mn3\mn\cdots\mnk}(x)=A^1_\pk(x).$$
\end{theorem}

\subsection{A classical pattern $\tau=2\mn1\mn3\mn\cdots\mnk$}
Using the same arguments as in the proof of Theorem~\ref{thbb1},
we get
\begin{theorem}\label{thbb4} For all $k\geq 3$,
\begin{enumerate}
\item $UD^1_{2\mn1\mn3\mn\cdots\mnk}(x)=\dfrac{1}{U_{k-1}^2\ttx}\left(\sum\limits_{j=0}^{k-3}U_{j+1}\ttx U_j\ttx -x\right)$;

\item $UU^1_{2\mn1\mn3\mn\cdots\mnk}(x)=\dfrac{1}{xU_{k-1}^2\ttx}\left(xU_{k-4}\ttx U_{k-2}\ttx +\sum\limits_{j=0}^{k-4} U_j\ttx
U_{j+1}\right)$.

\end{enumerate}
\end{theorem}

\section{Containing $1\mn3\mn2$ exactly once and another arbitrary
pattern}\label{sec5} In this section we consider those alternating
(up-up, up-down, down-up, down-down) permutations in $\SS_n$ that
contain $1\mn3\mn2$ exactly once and contain an arbitrary
generalized pattern $\tau$. We begin by setting some notation. Let
$A^1_{\tau;r}(x)$ be the generating function for the number of
alternating permutations in $A_n$ that contain $1\mn3\mn2$ exactly
once and contain $\tau$ exactly $r$ times. Moreover, We denote the
generating function for number of up-up (resp. up-down, down-up,
down-down) permutations in $\SS_n$ that contain $1\mn3\mn2$
exactly once and contain $\tau$ exactly $r$ times by
$UD^1_{\tau;r}(x)$ (resp. $UU^1_{\tau;r}(x)$, $DU^1_{\tau;r}(x)$,
$DD^1_{\tau;r}(x)$).

One can try to obtain results similar to Theorems~\ref{thaa1},
\ref{thc1}, and~\ref{thba1}, but expressions involved become
extremely cumbersome. So we just present the case of the up-down
permutations in $\SS_n$ that contain $1\mn3\mn2$ exactly once and
contain $\pk$ exactly once. Using Proposition~\ref{cv1} together
with Theorems~\ref{thaa1}, \ref{thc1}, and \ref{thbb1}, we get the
following result.

\begin{theorem}
Let $k\geq 2$, the generating function $UD_{\pk;1}^1(x)$ is given
by
$$\dfrac{1}{U_{k-1}^2\ttx}.\sum_{j=0}^{k-3}
\left[\dfrac{U_{j+1}\ttx\left( U_{j+1}\ttx+xU_j\ttx \right)
+2x\sum_{i=0}^{j-1} U_i\ttx U_{i+1}\ttx }{U_{j+1}\ttx U_{j+2}\ttx
} \right].$$
\end{theorem}

As a remark, the generating functions $UU_{\pk;1}^1(x)$,
$DU_{\pk;1}^1(x)$, and $DD_{\pk;1}^1(x)$ have more complicated
expressions since the equations, which these generating function
satisfied to, expressed in terms of $UD_{1\mn2\mn\cdots\mn
d;1}^1(x)$. For example, $UU_{\pk;1}^1(x)$ satisfies
$$\begin{array}{l}
UU_{\pk;1}^1(x)=xUD_{\pkk;1}^1(x)(1+UU_{\pk}(x))+xUD_{\pk}^1(x)UU_{\pk;1}(x)+\\
\qquad\qquad+xUD_{\pkk;1}(x)UU_{\pk}^1(x)+x(x+UD_{\pkk}(x))UU_{\pk;1}^1(x)+\\
\qquad\qquad\qquad+x^3UU_{\pkk;1}(x)(x+UD_{\pkk}(x))(1+UU_\pk(x))+\\
\qquad\qquad\qquad\ \ +x^3(1+UU_{\pkk}(x))UD_{\pkk;1}(x)(1+UU_\pk(x))+\\
\qquad\qquad\qquad\ \ \ \
+x^3(1+UU_{\pkk}(x))(x+UD_{\pkk}(x))UU_{\pk;1}(x).
\end{array}$$
\section{Further results}\label{sec6}
Here we present two directions to generalize and to extend the
results of the previous subsections.

\subsection{Statistics on alternating permutation that avoid $1\mn3\mn2$}
The first direction is to find a statistics on alternating
permutations in $A_n$ that avoid $1\mn3\mn2$.

\subsubsection{A classical pattern $\pk$} We define
$$\begin{array}{l}
UD_1(x_1,x_2,\dots)=\sum\limits_{n\geq0}\;\sum\limits_{\pi\in
UD_n(1\mn3\mn2)}\prod_{j\geq1}x_j^{1\mn2\mn\cdots\mn j(\pi)},\\
DU_1(x_1,x_2,\dots)=\sum\limits_{n\geq0}\;\sum\limits_{\pi\in
DU_n(1\mn3\mn2)}\prod_{j\geq1}x_j^{1\mn2\mn\cdots\mn j(\pi)},\\
UU_1(x_1,x_2,\dots)=\sum\limits_{n\geq0}\;\sum\limits_{\pi\in
UU_n(1\mn3\mn2)}\prod_{j\geq1}x_j^{1\mn2\mn\cdots\mn j(\pi)},\\
DD_1(x_1,x_2,\dots)=\sum\limits_{n\geq0}\;\sum\limits_{\pi\in
DD_n(1\mn3\mn2)}\prod_{j\geq1}x_j^{1\mn2\mn\cdots\mn j(\pi)},
\end{array}$$
where $1\mn2\mn\dots\mn j(\pi)$ is the number of occurrences of
$1\mn2\mn\cdots\mn j$ in $\pi$.

\begin{theorem}\label{st1} We have
\begin{itemize}
\item[(i)] The generating function $UD_1(x_1,x_2,\dots)$ is given
by
$$\frac{-\prod_{j\geq1}x_j^{2\binom{0}{j-1}}}
{\prod_{j\geq1}x_j^{\binom{0}{j-1}}-\dfrac{1}{\prod_{j\geq1}x_j^{\binom{0}{j-1}}-
\dfrac{\prod_{j\geq1}x_j^{2\binom{1}{j-1}}}{\prod_{j\geq1}x_j^{\binom{1}{j-1}}-
\dfrac{1}{\prod_{j\geq1}x_j^{\binom{1}{j-1}}-\dfrac{\prod_{j\geq1}x_j^{2\binom{2}{j-1}}}{\ddots}}}}}$$

\item[(ii)] The generating function $UU_1(x_1,x_2,\dots)$ is given
by
$$\frac{-1}{1-\dfrac{1}{x_1^2+x_1UD_1(x_1x_2,x_2x_3,\cdots)}}.$$

\item[(iii)] The generating function $DD_1(x_1,x_2,\dots)$ is given
by
$$\sum\limits_{n\geq1}\prod\limits_{d=1}^n
\dfrac{\prod_{j\geq1}x_j^{2\binom{d-1}{j-1}}}{1-\prod_{j\geq1}x_j^{2\binom{d-1}{j-1}}-\prod_{j\geq1}x_j^{\binom{d-1}{j-1}}
UD_1\left(\prod_{j\geq1}x_j^{\binom{d}{j}},\prod_{j\geq1}x_{j+1}^{\binom{d}{j}},\dots\right)}.$$

\item[(iv)] The generating function $DU_1(x_1,x_2,\cdots)$ is given
by
$$x_1UU_1(x_1,x_2,\cdots)+x_1DD_1(x_1x_2,x_2x_3,\cdots)(1+UU_1(x_1,x_2,\cdots)).$$
\end{itemize}
\end{theorem}
\begin{proof}
Proposition \ref{av1}(1) yields
$$UD_1(x_1,x_2,\cdots)=x_1(UD_1(x_1x_2,x_2x_3,\cdots)+x_1)(UD_1(x_1,x_2,\cdots)+x_1),$$
equivalently,
$$UD_1(x_1,x_2,\cdots)=\frac{-x_1^2}{x_1-\dfrac{1}{x_1+UD_1(x_1x_2,x_2x_3,\cdots)}}.$$
Hence, by induction, we get (i).

Proposition \ref{av1}(3) yields
$$UU_1(x_1,x_2,\cdots)=x_1(x_1+UD_1(x_1x_2,\cdots))(1+UU_1(x_1,x_2,\cdots)),$$
which equivalently to (ii).

Using Proposition \ref{av1}(4), we get
$$DD_1(x_1,x_2,\cdots)=x_1(x_1+UD_1(x_1,x_2,\cdots))(1+DD_1(x_1x_2,x_2x_3,\cdots)),$$
hence, by induction, we have (iii).

Finally, it is easy to check that Proposition \ref{av1}(2) yields
(iv).
\end{proof}

We denote the generating function
$C(x^2)-x=\frac{1-2x^2-\sqrt{1-4x^2}}{2x}$ by $\widehat{C}(x)$. As
an application to Theorem~\ref{st1}(i), we get the following:\\
(i) Let $x_1=x$ and $x_j=1$ for all $j\geq 2$. We have
$$UD_1(x,1,1,\dots)=\sum_{n\geq0}\sum_{\pi\in UD_n(1\mn3\mn2)}x^n=\widehat{C}(x),$$
as proved in Theorem~\ref{tha0}.\\
(ii) Let $\pi$ be any permutation. We say that $\pi_i$ is {\em
right to left maxima} if $\pi_i>\pi_j$ for all $j>i$. We denote
the number of right to left maxima of $\pi$ by $rlmax_\pi$. In
\cite{BCS}, there were proved that $rlmax_\pi=\sum\limits_{j\geq0}
(-1)^{j+1}1\mn2\mn\dots\mn j(\pi)$. Therefore
$$\sum_{n\geq0}\sum_{\pi\in
UD_n(1\mn3\mn2)}x^ny^{rlmax_\pi}=UD_1(xy,y^{-1},y,y^{-1},\dots)=\frac{x^2y^2(xy+\widehat{C}(x))}{1-x^2y^2-xy\widehat{C}(x)}.$$
In other words, the generating function for the number up-down
permutations in $UD_n(1\mn3\mn2)$ having $k$ right to left maxima
is given by
$$x^k\sum_{j=0}^{k-2}\binom{\frac{k-2+j}{2}}{\frac{k-2-j}{2}}\widehat{C}^j(x).$$
(iii) Let $\pi$ be any permutations. The number of increasing
subsequences in $\pi$ is given by
$inc_\pi=\sum_{j\geq0}1\mn2\mn\cdots\mn j(\pi)$. Hence,
$$\sum_{n\geq0}\sum_{\pi\in
UD_n(1\mn3\mn2)}x^ny^{inc_\pi}=Ud_1(xy,y,y,\ldots)=\dfrac{-x^2y^2}{xy-\dfrac{1}{xy-\dfrac{x^2y^4}{xy^2-\dfrac{1}{xy^2-\dfrac{x^2y^8}{\ddots}}}}}.$$
(iv) The generating function for the number of up-down
permutations that avoid $1\mn3\mn2$ and contain a prescribed
number of occurrences of the pattern $\pk$ is given by the
continued fraction in the statement of Theorem~\ref{st1}(i)
together with $x_1=x$, $x_k=y$, and $x_j=1$ for all $j\neq 1,k$.

As an application to Theorem~\ref{st1}(ii), we have the following:\\
(i) The generating function for the number of up-up permutations
in $\SS_n(1\mn3\mn2)$ is given by
$$UU_1(x,1,1,\dots)=\frac{1-\sqrt{1-4x^2}}{1+\sqrt{1-4x^2}}.$$
(ii) The distribution for the number of right to left maxima of
up-up permutations avoiding $1\mn3\mn2$ is given by
$$\begin{array}{l}
\sum\limits_{n\geq0}\sum\limits_{\pi\in
UU_n(1\mn3\mn2)}x^{|\pi|}y^{rlmax_\pi}=UU_1(xy,y^{-1},y,\dots)=\\
\qquad\qquad\qquad=\frac{-1}{1-\dfrac{1}{x^2y^2+xy\widehat{C}(x)}}=\sum\limits_{d\geq0}
\sum\limits_{j=0}^d\binom{d-j}{j}x^{2j}\widehat{C}^{d-2j}(x)y^d.
\end{array}$$
(iii) The generating function for the number of up-up permutations
that avoid $1\mn3\mn2$ and contain a prescribed number of
occurrences of the pattern $\pk$ is given by the statement of
Theorem~\ref{st1}(ii) together with $x_1=x$, $x_k=y$, and $x_j=1$
for all $j\neq 1,k$.

As an application to Theorem~\ref{st1}(iii), we have the following:\\
(i) The generating function for the number of down-down
permutations in $DD_n(1\mn3\mn2)$ is given by
$$DD_1(x,1,1,\dots)=\frac{1-\sqrt{1-4x^2}}{1+\sqrt{1-4x^2}}.$$
(ii) The distribution for the number of right to left maxima of
$\pi\in DD_n(1\mn3\mn2)$ is given by
$$\begin{array}{l}
\sum_{n\geq0}\sum\limits_{\pi\in
DD_n(1\mn3\mn2)}x^{|\pi|}y^{rlmax_\pi}
=DD_1(xy,y^{-1},y,\dots)=\dfrac{2x^2y^2}{\sqrt{1-4x^2}(1-xy\widehat{C}(x)-x^2y^2)}=\\
\qquad\qquad\qquad\qquad\qquad=\dfrac{2}{\sqrt{1-4x^2}}\sum\limits_{d\geq2}\sum\limits_{j=0}^{d-2}\binom{d-2-j}{j}\widehat{C}^{d-2j}(x)x^{d}y^{d}.
\end{array}$$
(iii) The generating function for the number of down-down
permutations avoiding $1\mn3\mn2$ and containing a prescribed
number of occurrences of the pattern $\pk$ is given by the
statement of Theorem~\ref{st1}(iii) together with $x_1=x$,
$x_k=y$, and $x_j=1$ for all $j\neq 1,k$.

As an application to Theorem~\ref{st1}(iv), we have the following:\\
(i) The generating function for the number of down-up permutations
in $DU_n(1\mn3\mn2)$ is given by
$$DU_1(x,1,1,\dots)=x\left( \frac{1-\sqrt{1-4x^2}}{2x^2}\right)^2-x.$$
(ii) The distribution for the number of right to left maxima of
$\pi\in DU_n(1\mn3\mn2)$ is given by
$$\begin{array}{l}
\sum\limits_{n\geq0}\sum\limits_{\pi\in
DU_n(1\mn3\mn2)}x^{|\pi|}y^{rlmax_\pi} =DU_1(xy,y^{-1},y,\dots)=\\
\qquad\qquad=\dfrac{1-\sqrt{1-4x^2}}{1+\sqrt{1-4x^2}}xy+\dfrac{2}{\sqrt{1-4x^2}}\sum\limits_{d\geq2}
x^dy^d\sum\limits_{j=0}^{d-1}\binom{d-1-j}{j}\widehat{C}^{d-1-2j}(x).
\end{array}$$
(iii) The generating function for the number of down-up
permutations avoiding $1\mn3\mn2$ and containing a prescribed
number of occurrences of the pattern $\pk$ is given by the
statement of Theorem~\ref{st1}(iv) together with $x_1=x$, $x_k=y$,
and $x_j=1$ for all $j\neq 1,k$.

\subsubsection{A generalized pattern $12\mn3\mn\cdots\mn k$} First of all, let us define
$$\begin{array}{l}
UD_2(x_1,x_2,\dots)=\sum\limits_{n\geq0}\;\sum\limits_{\pi\in
UD_n(1\mn3\mn2)}x_1^{|\pi|}\prod_{j\geq2}x_j^{12\mn3\mn\cdots\mn j(\pi)},\\
UU_2(x_1,x_2,\dots)=\sum\limits_{n\geq0}\;\sum\limits_{\pi\in
UU_n(1\mn3\mn2)}x_1^{|\pi|}\prod_{j\geq2}x_j^{12\mn3\mn\cdots\mn
j(\pi)},
\end{array}$$
where $12\mn3\mn\cdots\mn j(\pi)$ is the number of occurrences of
the generalized pattern $12\mn3\mn\dots\mn j$ in $\pi$. Similarly
to Theorem~\ref{st1} we get the following result.

\begin{theorem}\label{st2} We have
\begin{itemize}
\item[$(i)$] The generating function $UD_2(x_1,x_2,\dots)$ is given
by
$$\dfrac{-x_1^2\prod_{j\geq2}x_j^{\binom{0}{j-2}}}{x_1\prod_{j\geq2}x_j^{\binom{0}{j-2}}-
\dfrac{1}{x_1-\dfrac{x_1^2\prod_{j\geq2}x_j^{\binom{1}{j-2}}}{x_1\prod_{j\geq2}x_j^{\binom{1}{j-2}}-
\dfrac{1}{x_1-\dfrac{x_1^2\prod_{j\geq2}x_j^{\binom{2}{j-2}}}{x_1\prod_{j\geq2}x_j^{\binom{2}{j-2}}-
\dfrac{1}{x_1-\ddots} } } } } };$$

\item[$(ii)$] The generating function $UU_2(x_1,x_2,\dots)$ is
given by
$$\dfrac{-1}{1-\dfrac{1}{x_1^2x_2+x_1x_2UD_2(x_1,x_2x_3,x_3x_4,\dots)} };$$
\end{itemize}
\end{theorem}
As an application to Theorem \ref{st2}(i), we get that the number
of up-down permutations avoiding $1\mn3\mn2$ of length $2n+1$ with
$k$ rises is given by $C_n\delta_{k,(n-1)/2}$, where
$\delta_{a,b}=1$ if $a=b$,
otherwise $0$. \\
As an another application, we get that the generating function for
the number of up-down permutations avoiding $1\mn3\mn2$ and
containing a prescribed number of occurrences of the pattern
$12\mn3\mn4\mn\cdots\mnk$ is given by the statement of
Theorem~\ref{st2}(i) together with $x_1=x$, $x_k=y$, and $x_j=1$
for all $j\neq 1,k$. For example, for $k=3$ we get

\begin{equation}
\dfrac{-x^2}{x-\dfrac{1}{x-\dfrac{x^2y}
{xy-\dfrac{1}{x-\dfrac{x^2y^2}{xy^2-\ddots} } } } }.
\label{equd12x3}
\end{equation}

As an application to Theorem~\ref{st2}(ii), we get that the
generating function for the number of up-up permutations avoiding
$1\mn3\mn2$ and containing a prescribed number of occurrences of
the pattern $12\mn3\mn4\mn\cdots\mnk$ is given by the statement of
Theorem~\ref{st2}(ii) together with $x_1=x$, $x_k=y$, and $x_j=1$
for all $j\neq 1,k$. For example, for $k=3$ we get

\begin{equation}
\dfrac{-1}{1-\dfrac{1}{x^2-x\cdot\dfrac{x^2y}{xy-\dfrac{1}{x-\dfrac{x^2y^2}{xy^2-\ddots}}}}}.
\label{equu12x3}
\end{equation}

As a corollary to Theorem \ref{st2}(i,ii), we have the following
result.

\begin{corollary}
The generating function $\sum_{n\geq0}\sum_{\pi\in
A_n(1\mn3\mn2)}x_1^{|\pi|}\prod_{j\geq2}x_j^{12\mn3\mn\cdots\mn
j(\pi)}$ is given by
$$1+x_1+(1+x_1)UU_2(x_1,x_2,x_3,\dots).$$
\end{corollary}
\begin{proof}
Using Theorem~\ref{st2}(i,ii) we get that
$UD_2(x_1,x_2,x_3,\cdots)=x_1UU_2(x_1,x_2,x_3,\cdots)$. On the
other hand, by Identity~\ref{alt1} we have
$$\sum_{n\geq0}\sum_{\pi\in A_n(1\mn3\mn2)}x_1^{|\pi|}\prod_{j\geq2}x_j^{12\mn3\mn\cdots\mn j}=1+x_1+UD_2(x_1,x_2,x_3,\cdots)+UU_2(x_1,x_2,x_3,\cdots).$$
Hence, by combining the two equations above we get the desired
result.
\end{proof}

\subsubsection{A generalized pattern $21\mn3\mn\cdots\mn k$} First
of all, let us define another two objects
$$\begin{array}{l}
UD_3(x_1,x_2,\dots)=\sum\limits_{n\geq0}\sum\limits_{\pi\in
UD_n(1\mn3\mn2)}x_1^{|\pi|}\prod_{j\geq2}x_j^{21\mn3\mn\cdots\mn j(\pi)},\\
UU_3(x_1,x_2,\dots)=\sum\limits_{n\geq0}\sum\limits_{\pi\in
UU_n(1\mn3\mn2)}x_1^{|\pi|}\prod_{j\geq2}x_j^{21\mn3\mn\cdots\mn
j(\pi)},
\end{array}$$
where $21\mn3\mn\cdots\mn j(\pi)$ is the number of occurrences of
the generalized pattern $21\mn3\mn\dots\mn j$ in $\pi$. Similarly
to Theorem~\ref{st1} we get

\begin{theorem}\label{st3} We have
\begin{itemize}
\item[$(i)$]    $UD_3(x_1,x_2,\dots)=UD_2(x_1,x_2,\dots)$;

\item[$(ii)$]   $UU_3(x_1,x_2,\dots)=\frac{1}{x_2}UU_2(x_1,x_2,\dots)$.
\end{itemize}
\end{theorem}

Theorem~\ref{st3} and Identity~\ref{alt1} together with the fact
that $UD_2(x_1,x_2,x_3,\cdots)=x_1UU_2(x_1,x_2,x_3,\cdots)$ yield
the following corollary.

\begin{corollary}
The generating function $\sum_{n\geq0}\;\sum_{\pi\in
A_n(1\mn3\mn2)}x_1^{|\pi|}\prod_{j\geq2}x_j^{21\mn3\mn\cdots\mn
j(\pi)}$ is given by
$$1+x_1+\left(x_1+\frac{1}{x_2}\right)UU_2(x_1,x_2,x_3,\dots).$$
\end{corollary}
\subsection{Counting occurrences of $1\mn3\mn2$ in an alternating
permutation} The second direction is to find the generating
function $A_{\tau}(x;r)$ for the number of alternating permutation
in $\SS_n$ that contain $\tau$ exactly $r$ times. Moreover, to
find the generating functions $UU_{\tau}(x;r)$, $UD_{\tau}(x;r)$,
$DU_{\tau}(x;r)$, and $DD_{\tau}(x;r)$ for the number of up-up,
up-down, down-up, and down-down permutation in $\SS_n$ that
contain $\tau$ exactly $r$ times, respectively.

For $\tau=1\mn3\mn2$, the block decomposition approach~\cite{MVr}
in the alternating (up-up, up-down, down-up, down-down)
permutations gives a complete answer for any given $r$ (as
described in Section~2 for $r=0$, and Section~4 for $r=1$), and we
get the following result.

\begin{theorem}\label{rtimes} We have
$$\begin{array}{ll}
(1)& \left\{\begin{array}{l}
        UD_{1\mn3\mn2}(x;0)=\frac{1-2x^2}{2x}-\frac{1}{2x}(1-4x^2)^{\frac{1}{2}};\\
        \\
        UD_{1\mn3\mn2}(x;1)=\frac{-1}{2x}+\frac{1-2x^2}{2x}(1-4x^2)^{\frac{-1}{2}};\\
        \\
        UD_{1\mn3\mn2}(x;2)=\frac{1}{2x}-\frac{1-6x^2+6x^4}{2x}(1-4x^2)^{\frac{-3}{2}};\\
        \\
        UD_{1\mn3\mn2}(x;3)=\frac{-2(1-x^2)}{x}+\frac{2-22x^2+80x^4-98x^6+16x^8}{x}(1-4x^2)^{\frac{-5}{2}},
        \end{array}\right.\\

&\\
(2)& \left\{\begin{array}{l}
        UU_{1\mn3\mn2}(x;0)=\frac{1-2x^2}{2x^2}-\frac{1}{2x^2}(1-4x^2)^{\frac{1}{2}};\\
        \\
        UU_{1\mn3\mn2}(x;1)=\frac{x^2-1}{x^2}+\frac{1-3x^2}{x^2}(1-4x^2)^{\frac{-1}{2}};\\
        \\
        UU_{1\mn3\mn2}(x;2)=\frac{4-5x^2}{2x^2}-\frac{4-29x^2+54x^4-16x^6}{2x^2}(1-4x^2)^{\frac{-3}{2}};\\
        \\
        UU_{1\mn3\mn2}(x;3)=\frac{13-11x^2+2x^4}{2x^2}+\frac{13-152x^2+612x^4-940x^6+384x^8}{2x^2}(1-4x^2)^{\frac{-5}{2}},
        \end{array}\right.\\
&\\
(3) & \left\{\begin{array}{l}
        \\
        A_{1\mn3\mn2}(x;0)=\frac{1+x-2x^2-2x^3}{2x^2}-\frac{1+x}{2x^2}(1-4x^2)^{\frac{1}{2}};\\
        \\
        A_{1\mn3\mn2}(x;1)=\frac{2x^2-x-2}{2x^2}+\frac{2+x-6x^2-2x^3}{2x^2}(1-4x^2)^{\frac{-1}{2}};\\
        \\
        A_{1\mn3\mn2}(x;2)=\frac{4+x-5x^2}{2x^2}-\frac{4+x-29x^2-6x^3+54x^4+6x^5-16x^6}{2x^2}(1-4x^2)^{\frac{-3}{2}};\\
        \\
        A_{1\mn3\mn2}(x;3)=\frac{13-4x-11x^2+4x^3+2x^4}{2x^2}+\\
        \quad\quad\qquad\qquad\qquad+\frac{13+4x-152x^2-44x^3+612x^4+160x^5-940x^6-196x^7+384x^8+32x^9}{2x^2}(1-4x^2)^{\frac{-5}{2}}.
        \end{array}\right.
\end{array}$$
\end{theorem}

Finally, we conclude with some open problems suggested by the
results in the previous sections. For example: (1) We found that
the number of alternating permutations in $A_n(1\mn3\mn2, \pk)$ is
equal to the number of alternating permutations in $A_n(1\mn3\mn2,
2\mn1\mn3\mn\cdots\mnk)$ (see Corollary~\ref{ccaa1}). The
question, if there is a bijective combinatorial proof for this
property. (2) If there exists an explicit formula for
$A_{1\mn3\mn2}(x;r)$ for any given $r\geq0$. (3) If there exist
bijective combinatorial proofs for the formulas in the statement
of Theorem~\ref{rtimes}.


\begin{thebibliography}{WWW}
\bibitem[A1]{AP1}
D.~Andr\'e, Developments de $\sec x$ et $\tan x$, {\em C.R. Acad.
Sci.}, Paris {\bf 88} (1879), 965--967.

\bibitem[A2]{AP2}
D.~Andr\'s, memire sur les permutations altern\'ees, {\em J.
Math.} {\bf 7} (1881), 167--184.

\bibitem[Bo]{B1}
M.~B\'ona, The permutation classes equinumerous to the smooth
class, {\em Electron. J. Combin.} {\bf 5} (1998), \#R31.

\bibitem[BS]{BS}
E.~Babson and E.~Steingr\'\i msson, Generalized permutation
patterns and a classification of the Mahonian statistics, {\em
S\'eminaire Lotharingien de Combinatoire}, B44b:18pp, (2000).

\bibitem[BCS]{BCS}
P.~Br\"and\'en, A.~Claesson, and E.~Steingr\'\i msson, Continued
fractions and increasing subsequences in permutations, {\em Discr.
Math.}, to appear.

\bibitem[C]{C}
A.~Claesson, Generalised pattern avoidance, {\em Europ. J.
Combin.} {\bf 22} (2001), 961--973.

\bibitem[CW]{CW}
T.~Chow and J.~West, Forbidden subsequences and Chebyshev
polynomials, {\em Discr. Math.} {\bf 204} (1999), 119--128.

\bibitem[Kn]{Kn}
D.E.~Knuth, The Art of Computer Programming, 2nd ed. Addison
Wesley, Reading, MA (1973).

\bibitem[Km]{Km}
D.~Kremer, Permutations with forbidden subsequences and a
generalized Schr\"oder number, {\em Discr. Math.} {\bf 218}
(2000), 121--130.

\bibitem[Kr]{Kr}
C.~Krattenthaler, Permutations with restricted patterns and Dyck
paths, {\em Adv. Appl. Math.} {\bf 27} (2001), 510--530.

\bibitem[M1]{Mg1}
T.~Mansour, Continued fractions and generalized patterns, {\em
Europ. J. Combin.} {\bf 23:3} (2002), 329--344.

\bibitem[M2]{Mg2}
T.~Mansour, Continued fractions, statistics, and generalized
patterns, {\em Ars Combinatorica}, to appear (2002), preprint
CO/0110040.

\bibitem[M3]{Mg3}
T.~Mansour, Restricted $1$-$3$-$2$ permutations and generalized
patterns, {\em Annals of Combinatorics} {\bf 6} (2002), 65--76.

\bibitem[MV1]{MV1}
T.~Mansour and A.~Vainshtein, Restricted permutations, continued
fractions, and Chebyshev polynomials {\em Electron. J. Combin.}
\textbf{7} (2000), \#R17.

\bibitem[MV2]{MV2}
T.~Mansour and A.~Vainshtein, Restricted 132-avoiding
permutations, {\em Adv. Appl. Math.} {\bf 126} (2001), 258--269.

\bibitem[MV3]{MV3}
T.~Mansour and A.~Vainshtein, Layered restrictions and Chebychev
polynomials, {\em Annals of Combin.} {\bf 5} (2001), 451--458.

\bibitem[MV4]{MV4}
T.~Mansour and A.~Vainshtein, Restricted permutations and
Chebyshev polynomials, {\em S\'eminaire Lotharingien de
Combinatoire} {\bf 47} (2002), Article B47c.

\bibitem[MV5]{MVr}
Counting occurrences of $132$ in a permutation, {\em Adv. Appl.
Math.} {\bf 28:2} (2002), 185--195.

\bibitem[R]{R}
A.~Robertson, Permutations containing and avoiding 123 and 132
patterns, {\em Disc. Math. and Theo. Comp. Sci.} {\bf 3} (1999),
151--154.

\bibitem[Ri]{Ri}
Th.~Rivlin, Chebyshev polynomials. From approximation theory to
algebra and number theory, John Wiley, New York (1990).

\bibitem[RWZ]{RWZ}
A.~Robertson, H.~Wilf, and D.~Zeilberger, Permutation patterns and
continuous fractions, {\em Elec. J. Combin.} {\bf 6} (1999),
\#R38.

\bibitem[SP]{SP}
N.J.A.~Sloane and S.~Plouffe, {\em The Encyclopedia of Integer
Sequences}, Academic Press, New York (1995).

\bibitem[SS]{SS}
R.~Simion, F.W.~Schmidt, Restricted Permutations, {\em Europ. J.
Combin.} {\bf 6} (1985), 383--406.

\bibitem[W]{W}
J.~West, Generating trees and forbidden subsequences, {\em Discr.
Math.} {\bf 157} (1996), 363--372.
\end{thebibliography}
\end{document}